\documentclass[12pt]{article}
\usepackage{amsmath,amssymb,amsthm}
\usepackage{hyperref}
\usepackage{datetime}
\usepackage{units}
\usepackage{color}
\usepackage[T1]{fontenc}
\usepackage[utf8]{inputenc}
\usepackage{authblk}
\usepackage{bm,latexsym,mathrsfs,enumerate}
\usepackage{color}

\title{A New Binary BBP-type Formula for $\sqrt 5\,\log\phi$ \thanks{%
MSC 2010: 11Y60, 30B99}}
\author[]{Kunle Adegoke\thanks{adegoke00@gmail.com}}
\affil{Department of Physics and Engineering Physics, \mbox{Obafemi Awolowo University, Ile-Ife, 220005 Nigeria}}

\theoremstyle{plain}
\numberwithin{equation}{section}

\begin{document}
\date{}
\maketitle
\begin{abstract}
\noindent Hitherto only a base 5 BBP-type formula is known for $\sqrt 5\log\phi$, where \mbox{$\phi=(\sqrt 5+1)/2$}, the golden ratio, ( i.e. Formula 83 of the April 2013 edition of Bailey's Compendium of \mbox{BBP-type} formulas). In this paper we derive a new binary BBP-type formula for this constant. The formula is obtained as a particular case of a BBP-type formula for a family of logarithms.
\end{abstract}
\tableofcontents

\section{Introduction}

A Bailey-Borwein-Ploufe (BBP) type formula for a mathematical constant $c$ has the form

\[
c= \sum\limits_{k = 0}^\infty  {\frac{1}{{b^k }}\sum\limits_{j = 1}^l {\frac{{a_j }}{{(kl + j)^s }}} }\,,
\]

where $s$, $b$ and $l$ are integers: the degree, base and length of the formula, respectively and the $a_j$ are rational numbers.

\bigskip

Such a formula allows the extraction of the individual base $b$ digits of a mathematical constant without the need to compute the previous digits. The original BBP formula, discovered in 1996~\cite{bbp97}, allows the extraction of the binary or hexadecimal digits of the constant $\pi$. Many such formulas have since been discovered and can be found in Bailey's Online Compendium of BBP-type formulas~\cite{bailey13} and in the references therein. Another online Compendium is also being maintained by the CARMA Institute~\cite{carma}. As at the time of writing this paper only a base $5$ BBP-type formula is known for the mathematical constant $\sqrt 5\log\phi$, where \mbox{$\phi=(\sqrt 5+1)/2$} is the golden ratio. This formula is listed as Formula~83 in the current edition of Bailey's Compendium.

\bigskip

In this present paper we derive a new binary BBP-type formula for $\sqrt 5\log\phi$. The formula is presented as a particular case of a more general BBP-type formula derived for a family of logarithms.

\section{Notation}
 
The first degree polylogarithm function, which we employ in this paper, is defined by

\[
{\rm Li}_1 [z] = \sum\limits_{r = 1}^\infty  {\frac{{z^r }}{{r }}}=-\log(1-z),\quad |z|\le 1,\;z\ne 1\,. 
\]

For $q,x\in\mathbb{R}$, we have the identities
\begin{equation}\label{equ.ne7tbp0}
\arctan \left( {\frac{{q\sin x}}{{1 - q\cos x}}} \right) = {\rm Im\;Li}_1 \left[ {q\exp (ix)} \right]=\sum_{r = 1}^\infty { {\frac{{q^r \sin rx}}{r}}}
\end{equation}

and

\begin{equation}\label{equ.uumx4p6}
- \frac{1}{2}\log \left( {1 - 2q\cos x + q^2 } \right) = {\rm Re\;Li}_1 \left[ {q\exp (ix)} \right]=\sum_{k = 1}^\infty {{\frac{{q^r \cos rx}}{r}}}\,.
\end{equation}

\bigskip

The BBP-type formulas in this paper will be given in the standard notation for BBP-type formulas, introduced in~\cite{bailey01}:

\[
\sum\limits_{k = 0}^\infty  {\frac{1}{{b^k }}\sum\limits_{j = 1}^l {\frac{{a_j }}{{(kl + j)^s }}} }\equiv P(s,b,l,A)\,,
\]

where $s$, $b$ and $l$ are integers, and \mbox{$A = (a_1, a_2,\ldots, a_l)$} is a vector of integers.

\section{A general BBP-type Formula for a family of logarithms}

\subsection*{Theorem}
{\em For any nonzero integer $t$, the following BBP-type formula holds}

\begin{equation}\label{equ.iypnkgl}
\begin{split}
&\sqrt 5\,\tanh^{-1}\left\{\left(\frac{1-t+2t^2}{1-t+3t^2-2t^3+4t^4}\right)t\sqrt 5\right\}\\
&=\frac{5}{2^{20}t^{39}}P(1,2^{20} t^{40},40,(2^{19}\, t^{38},0,2^{18}\, t^{36},2^{18}\, t^{35},0,0,-2^{16}\, t^{32},2^{16}\, t^{31},2^{15}\, t^{30},\\ &\qquad\quad 0,-2^{14}\, t^{28},-2^{14}\, t^{27},2^{13}\, t^{26},0,0,-2^{12}\, t^{23},-2^{11}\, t^{22},0,-2^{10}\, t^{20},0,-2^9\, t^{18},\\ &\qquad\quad 0,-2^8\, t^{16},-2^8\, t^{15},0,0,2^6\, t^{12},-2^6\, t^{11},-2^5\, t^{10},0,2^4\, t^{8},2^4\, t^{7},-2^3\, t^{6},0,0,\\ &\qquad\qquad 2^2\, t^{3},2\, t^{2},0,1,0))\,.
\end{split}
\end{equation}



\begin{proof}
Consider the following identity which holds for $t\in\mathbb{R},\,t\ne 0$:
\begin{equation}\label{equ.dksv20r}
\begin{split}
&\tanh ^{ - 1} \left\{ {\left( {\frac{{1 - t + 2t^2 }}{{1 - t + 3t^2  - 2t^3  + 4t^4 }}} \right)t\sqrt 5 } \right\}\\
&\quad ={\rm Re\,Li}_1 \left[ {\frac{1}{{t\sqrt 2 }}\exp \left( {\frac{{i\pi }}{{20}}} \right)} \right] - {\rm Re\,Li}_1 \left[ {\frac{1}{{t\sqrt 2 }}\exp \left( {\frac{{7i\pi }}{{20}}} \right)} \right]\\
&\qquad + {\rm Re\,Li}_1 \left[ {\frac{1}{{t\sqrt 2 }}\exp \left( {\frac{{9i\pi }}{{20}}} \right)} \right] - {\rm Re\,Li}_1 \left[ {\frac{1}{{t\sqrt 2 }}\exp \left( {\frac{{17i\pi }}{{20}}} \right)} \right]\,.
\end{split}
\end{equation}

It is straightforward to verify equation~\eqref{equ.dksv20r} by the use of the first equality in equation~\eqref{equ.uumx4p6}. Using the second equality in equation~\eqref{equ.uumx4p6} and trigonometric addition rules, equation~\eqref{equ.dksv20r} can also be written as

\begin{equation}\label{equ.wp3hkgy}
\begin{split}
\tanh ^{ - 1} \left\{ {\left( {\frac{{1 - t + 2t^2 }}{{1 - t + 3t^2  - 2t^3  + 4t^4 }}} \right)t\sqrt 5} \right\} &= \sum\limits_{r = 1}^\infty  {\frac{4}{{r\sqrt {2^r t^{2r} } }}\sin \left( {\frac{{r\pi }}{5}} \right)\sin \left( {\frac{{2r\pi }}{5}} \right)\cos \left( {\frac{{r\pi }}{4}} \right)}\\
&= \sum\limits_{r = 1}^\infty  {\frac{1}{{r\sqrt {2^r t^{2r} } }}f(r)}\,, 
\end{split}
\end{equation}
where we have defined a periodic function

\[
f(r)={4\sin \left( {\frac{{r\pi }}{5}} \right)\sin \left( {\frac{{2r\pi }}{5}} \right)\cos \left( {\frac{{r\pi }}{4}} \right)},\;r\in\mathbb{Z}\,.
\]

Since $f(r)\in \left\{0,\pm \sqrt 5/\sqrt 2,\pm\sqrt 5\right\}$ and since $f(40k+j)=f(j)$ for integers $k$ and $j$, we can convert the above single sum to an equivalent double sum by setting $r=40k+j$ in equation~\eqref{equ.wp3hkgy}, obtaining

\begin{equation}\label{equ.kc2o8gm}
\begin{split}
\tanh ^{ - 1} \left\{ {\left( {\frac{{1 - t + 2t^2 }}{{1 - t + 3t^2  - 2t^3  + 4t^4 }}} \right)t\sqrt 5 } \right\} &= \sum\limits_{k = 0}^\infty  {\left\{ {\frac{1}{({2^{20} t^{40})^k }}\sum\limits_{j = 1}^{40} {\left( {\frac{{f(j)}}{{t^j }}\frac{1}{{\sqrt {2^j } }}\frac{1}{{40k + j}}} \right)} } \right\}}\\
&= \frac{{\sqrt 5 }}{{2^{20} t^{39} }}\sum\limits_{k = 0}^\infty  {\left\{ {\frac{1}{{(2^{20} t^{40})^k }}\sum\limits_{j = 1}^{40} {\left( {\frac{{a_j }}{{40k + j}}} \right)} } \right\}}\,,
\end{split} 
\end{equation}

where the integers (for $t\in\mathbb{Z},t\ne 0$) $a_j$, $j=1,2,\ldots, 40$ are given by

\[
a_j  = \frac{{f(j)}}{5}t^{39 - j} \sqrt 5 \sqrt {2^{40 - j} }\,.
\]

Putting the explicit values of $a_j$ into equation~\eqref{equ.kc2o8gm} the theorem is proved.
\end{proof}

\subsection*{Remark}



The theorem (equation~\eqref{equ.iypnkgl}) is actually true for any nonzero complex number $t$, as is readily established by the Principle of Permanence for analytic functions~\cite{asimov}. Technically speaking, however, equation~\eqref{equ.iypnkgl} is BBP-type only if $t$ is a nonzero integer.

\section{A binary BBP-type Formula for $\sqrt 5\,\log\phi$}

Setting $t=1$ in equation~\eqref{equ.iypnkgl} gives
\subsection*{Corollary}

\begin{equation}\label{equ.h46nul5}
\begin{split}
\sqrt 5\log\phi&=\frac{5}{3\cdot 2^{20}}P(1, 2^{20}, 40, (2^{19}, 0, 2^{18}, 2^{18}, 0, 0, -2^{16}, 2^{16}, 2^{15}, 0, -2^{14}, -2^{14},\\ &\qquad \quad 2^{13}, 0, 0, -2^{12}, -2^{11}, 0, -2^{10}, 0, -2^{9}, 0, -2^{8}, -2^{8}, 0, 0, 2^{6}, -2^{6},\\ &\qquad -2^{5}, 0, 2^{4}, 2^{4}, -2^{3}, 0, 0, 2^{2}, 2, 0, 1, 0))\,.
\end{split}
\end{equation}


\section{Conclusion}

We have derived a new BBP-type formula for a family of logarithms. A binary BBP-type formula for $\sqrt 5\log\phi$ is obtained as a particular case of the result.

\section*{Acknowledgments}

The author thanks the anonymous referee for an excellent review and for information on the Principle of Permanence for analytic functions.

\bibliographystyle{unsrt}

\end{document}